\documentclass{conm-p-l}
\usepackage{mathrsfs}
\usepackage[percent]{overpic}

\theoremstyle{remark}

\numberwithin{equation}{section}
\numberwithin{figure}{section}

\newcommand\thmenv[3]{\begin{#1} \label{#1:#2} #3 \end{#1}}

\def\procl#1.#2 #3\endprocl{%
       \ifx#1t\thmenv{thm}{#2}{#3}\fi
       \ifx#1p\thmenv{pro}{#2}{#3}\fi
       \ifx#1c\thmenv{cor}{#2}{#3}\fi
       \ifx#1q\thmenv{question}{#2}{#3}\fi
       \ifx#1r\thmenv{remark}{#2}{{\rm #3}}\fi
    }%

\def\rref#1.#2/{%
      \ifx #1sSection~\ref{s.#2}\fi
      \ifx #1tTheorem~\ref{thm:#2}\fi  
      \ifx #1cCorollary~\ref{cor:#2}\fi 
      \ifx #1rRemark~\ref{remark:#2}\fi 
      \ifx #1fFigure~\ref{fig:#2}\fi
      \ifx #1e(\ref{e.#2})\fi
      \ifx #1b\cite{#2}\fi
      \ifx #1B\cite{#2}\fi
        }

\def\rlabel #1 #2{\begin{equation} \label{#1} #2 \end{equation}}

\newcommand\rproof{\begin{proof}}

\newcommand\eqaln[1]{\begin{align*} #1 \end{align*}}

\newcommand\textfrac{\frac}

\newcommand\Qed{\end{proof}}

\newcommand\bsection[2]{\bigbreak\section{#1}\label{#2}}

\newcommand\gh{G}  
\newcommand\verts{\mathsf{V}}
\newcommand\vertex{\mathsf{V}}
\newcommand\edges{\mathsf{E}}
\newcommand\Tr{\mathop{\rm Tr}}  
\newcommand\ip[1]{(\changecomma #1)}
\newcommand\bigip[1]{\bigl(\changecomma #1\bigr)}

\newcommand\bigpip[1]{\bigl(\changecomma #1\bigr)}

\def\changecomma#1,{#1,\,}
\def\bigchangecomma#1,{#1,\;}
\def\leftchangecomma#1,{#1,\ }

\newcommand\bp{o}
\newcommand\alg{\mathsf{Alg}}  
\newcommand\algA{\mathsf{A}}  
\newcommand\affalg{\mathsf{AffAlg}}  
\newcommand\LL{{\mathcal L}}  
\newcommand\etail[1]{#1^-}   
\newcommand\ehead[1]{#1^+}   
\newcommand\HilG{\mathscr{G}}  
\newcommand\HilH{\mathscr{H}}  
\newcommand\dint{\int^{\oplus}}
\newcommand\all[1]{\forall #1\enspace}
\newcommand\I[1]{{\bf 1}_{#1}}
\newcommand\II[1]{\I{\{#1\}}}
\newcommand\st{\, ; \;}  
\newcommand\Aut{{\rm Aut}} 

\DeclareMathOperator{\bP}{\mathbf{P}\mathopen{}}
\DeclareMathOperator{\E}{\mathbf{E}\mathopen{}}

\newcommand\Psub [1]{\bP_{\! #1}}

\newcommand\Psubbig [2]{\Psub {#1}\mkern-1.5mu\bigl[#2\bigr]}

\newcommand\PBig[1]{\bP\mkern-.5mu\Bigl[#1\Bigr]}

\newcommand\Ebig[1]{\E\mkern-1.5mu\bigl[#1\bigr]}
\newcommand\Esubbig [2]{\E_{#1}\mkern-1.5mu\bigl[#2\bigr]}

\newcommand\EBig[1]{\E\mkern-1.5mu\Bigl[#1\Bigr]}

\newcommand\dfn[1]{\textit{\textbf{#1}}}

\newcommand\Z{{\mathbb Z}}
\newcommand\R{{\mathbb R}}

\begin{document}

\title[Monotonicity of return probabilities]{Monotonicity of average return probabilities\\
for random walks in random environments}

\author{Russell Lyons}
\address{Department of Mathematics, Indiana
University, 831 E. 3rd St., Bloomington, IN 47405-7106}
\email{rdlyons@indiana.edu}
\thanks{Partially supported by the National
Science Foundation under grant DMS-1612363
and by Microsoft Research.}
\urladdr{http://pages.iu.edu/~rdlyons/}

\date{May 10, 2017.}

\subjclass[2010]{Primary 
 60K37, 
 60J35; 
Secondary
 05C80, 
 05C81
}
\keywords{Continuous time, Markov chains.}

\begin{abstract}
We extend a result of Lyons (2016) from fractional tiling of
finite graphs to a version for infinite random graphs.
The most general result is as follows.
Let $\bP$ be a unimodular probability measure on rooted networks
$(G, o)$ with positive weights $w_G$ on
its edges and with a percolation subgraph $H$ of $G$ with positive weights
$w_H$ on its edges. 
Let $\Psub {(G, o)}$ denote the conditional law of $H$ given $(G, o)$.
Assume that $\alpha := \Psubbig {(\gh, \bp)}{o \in \verts(H)} > 0$ is a
constant $\bP$-a.s.
We show that if $\bP$-a.s.\
whenever $e \in \edges(\gh)$ is adjacent to $\bp$,
\[
\Esubbig {(\gh, \bp)}{w_H(e) \bigm| e \in \edges(H)} 
\Psubbig {(\gh, \bp)}{e \in \edges(H) \bigm| \bp\in \verts(H)} 
\le
w_G(e)
\,,
\]
then 
\[
\forall t > 0 \quad \Ebig{p_t(o; G)} \le \Ebig{p_t(o; H) \bigm| o \in \verts(H)}
\,.
\]
\end{abstract}

\maketitle

\bsection{Introduction}{s.intro}

Associated to a graph with nonnegative numbers on its edges such that the
sum of numbers of edges incident to each given vertex is finite, there is a
continuous-time random walk that, when at a vertex $x$, crosses each edge
$e$ incident to $x$ at rate equal to the number on $e$. When all rates
equal 1, this is called continuous-time simple random walk. 
In general, the rate at which the random walk leaves $x$ equals the sum of
the numbers on the edges incident to $x$.

It is well known and easy to prove
that every such (weighted) random walk has the property
that the probability of return to the starting vertex is a decreasing
function of time.
Equivalently, the return probability at any fixed time decreases if all the
rates are increased by the same factor.
However, the return probability is {\it not\/} a decreasing function of the
set of rates in general. Indeed, the behavior of the return probabilities
is not intuitive; a small example is shown in \rref f.3pathex/.
Examples show that the return probability to a vertex $x$
need not be monotonic even when
rates are changed only on edges not incident to $x$.
On the other hand, on a finite graph,
the {\it average} of the return probabilities {\it is}
decreasing in the rates,
as shown by Benjamini and Schramm (see Theorem 3.1 of \rref
b.HeicklenHoffman/).
Recall that on a finite graph, the stationary measure for this continuous-time
random walk is uniform on the vertices.

\begin{figure}
\begin{overpic}[width=4.5truein,tics=5]{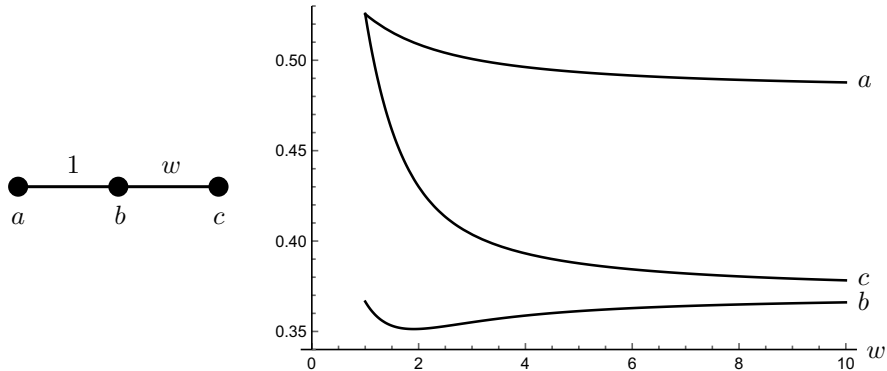}
   \put (8,23) {$1$}
   \put (19,23) {$w$}
   \put (1.5,17) {$a$}
   \put (13.5,17) {$b$}
   \put (25,17) {$c$}
   \put (101,1.5) {$w$}
   \put (100,33) {$a$}
   \put (100,7) {$b$}
   \put (100,10) {$c$}
\end{overpic}
\caption{
The return probabilities {\it at time 1}
of vertices $a$, $b$, and $c$ in 
the graph on the left as the rate $w$ varies.
}
\label{fig:3pathex}  
\end{figure}

In Theorems 4.1 and 4.2 of
\rref b.Lyons:treedom/, we extended and strengthened
the theorem of Benjamini and Schramm to the case of graphs of different
sizes and even to the case of one graph $G$ that is ``fractionally tiled"
by a set of subgraphs $H_i$ with a certain condition on the edge weights of $G$
and $H_i$. Our purpose here is to establish a version of those results for
infinite graphs.

For a very simple example of our results here,
consider the square lattice $\Z^2$ and the
subgraph $H$ formed by deleting every vertex both of whose coordinates are
odd; see \rref f.Z2andH/.
There are four subgraphs of $\Z^2$ that are isomorphic to $H$.
Considering those four copies of $H$, we find that each vertex of $\Z^2$ is
covered three times, once by a vertex of degree 4 and twice by a
vertex of degree 2. An appropriate average return probability in $H$
is thus $1/3$ that of a vertex of degree 4 plus $2/3$ that of a vertex of
degree 2.
Consider continuous-time simple random walk on each graph, where edges
are crossed at rate 1; the return probabilities are shown in \rref
f.rtnZ2eachH/.
As illustrated in \rref f.rtnZ2andH/, we have for all $t \ge 0$, 
\[
p_t\bigl((0, 0); \Z^2\bigr)
\le
\textfrac13 p_{3t/2}\bigl((0, 0); H\bigr) +
\textfrac23 p_{3t/2}\bigl((0, 1); H\bigr) 
\,.
\]
Effectively, we have used rates $3/2$ on every edge of $H$.
This inequality follows from \rref c.transitive/.
It is sharp in the following sense:
if $3t/2$ is replaced by $\beta t$ for some $\beta > 3/2$, then
the resulting inequality fails for some $t > 0$.

\begin{figure}
\includegraphics[width=2truein]{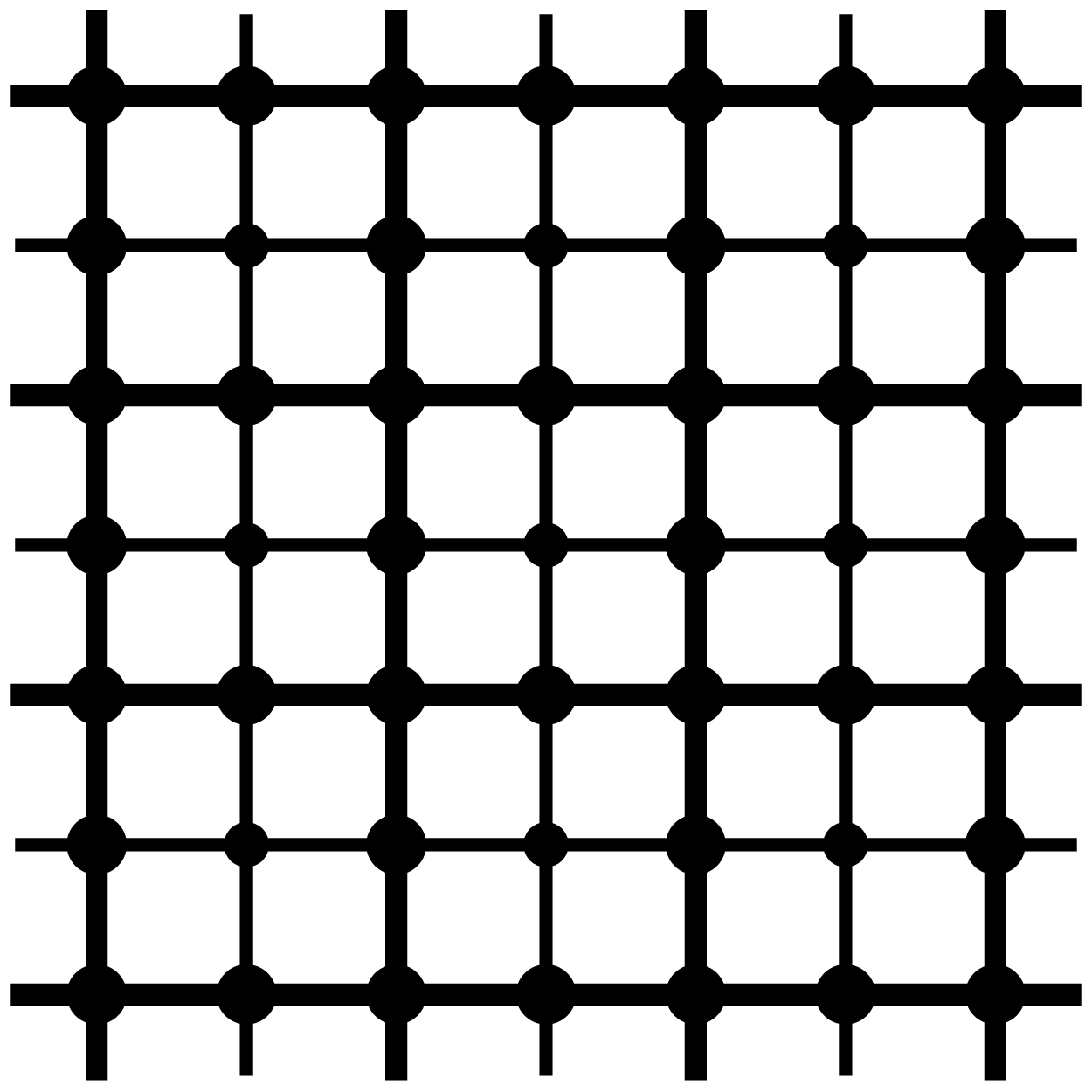}
\caption{The square lattice $\Z^2$ and the subgraph $H$, drawn thicker,
formed by deleting every vertex both of whose coordinates are odd.}
\label{fig:Z2andH}
\end{figure}

\vfill

\begin{figure}
\begin{overpic}[width=3truein]{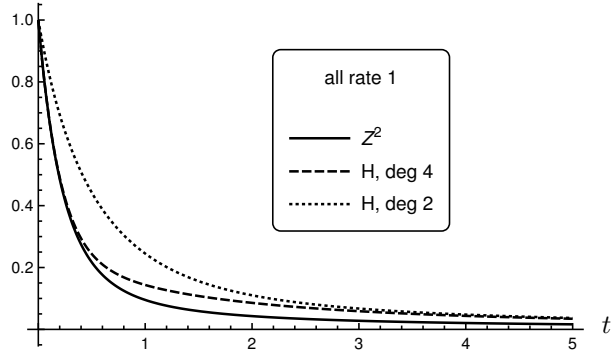}
\put (103,03) {$t$}
\end{overpic}
\caption{Return probabilities for
continuous-time simple random walk on each graph.}
\label{fig:rtnZ2eachH}
\end{figure}

\vfill

\begin{figure}
\begin{overpic}[width=3truein]{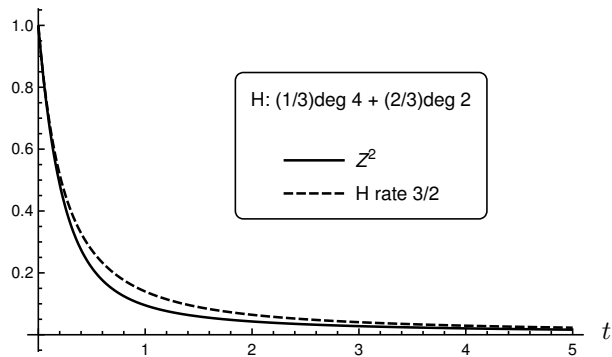}
\put (103,03) {$t$}
\end{overpic}
\caption{Comparison of continuous-time random walk with rates 1
on $\Z^2$ and rates $3/2$ on $H$, the latter
averaged over the starting point.}
\label{fig:rtnZ2andH}
\end{figure}


This particular example can be easily derived from Theorem 4.2 of \rref
b.Lyons:treedom/. With some more work, so can all the results here when the
unimodular
probability measures involved are sofic. However, the general case (which
is not known to be sofic) does not follow from earlier work. 
Nevertheless,
our results and proofs are modeled on Theorem 4.2 of \rref b.Lyons:treedom/.
The challenge here was to formulate the proper statements for
infinite graphs and to make the appropriate adjustments to the proofs
required for using direct integrals instead of direct sums.

For a more complicated example of our results,
suppose that $\gh$ is the usual nearest-neighbor graph on
$\Z^d$ ($d \ge 2$) and $H$ is the infinite cluster of supercritical
Bernoulli (site or bond)
percolation on $\gh$. Let $\delta := \Ebig{\deg_H(\bp) \bigm| \bp
\in H}/(2d) \in (0, 1)$. Then 
\[
\all {t \ge 0} \quad p_t(o; \Z^d) 
\le
\Ebig{p_{t/\delta}(o; H) \bigm| o \in \verts(H)}
\,.
\]
This is obtained by using $w_H \equiv 1/\delta$ in \rref c.transitive/.
The preceding inequality is false for any larger value of $\delta$.

\bsection{Statements of Results and Background}{s.ftr}

Let $G$ be a simple, locally finite graph with weights $w_G(e) \ge 0$
on the edges $e$.
Consider the continuous-time random walk on $G$ where edge $e$ is crossed
at rate $w_G(e)$ when the walk is incident to $e$.
Let $p_t(x; G)$ denote the probability that a random walk
started at $x$ is found at $x$ at time $t$.
If $\Delta_G$ is the corresponding Laplacian, i.e., $\Delta_G(x, y) :=
- w(e)$ when $e$ is an edge joining $x$ and $y$ with weight $w(e)$, all
other off-diagonal elements of $\Delta_G$ are 0, and the row sums are 0,
then $p_t(x; G)$ is the $(x, x)$-entry of $e^{-t \Delta_G}$.
If the entries of $\Delta_G$ are unbounded, then we take the minimal Markov
process, which dies after an explosion.
The infinitesimal generator is then the self-adjoint extension of
$-\Delta_G$ (for uniqueness of the extension, see \rref b.HKMW/).
For the definition of unimodular in our context, see \rref b.AL:urn/.

\procl t.avgrtn
Let $\bP$ be a unimodular probability measure on rooted networks
$(G, o)$ with positive weights $w_G$ on
its edges and with a percolation subgraph $H$ of $G$ with positive weights
$w_H$ on its edges. 
Let $\Psub {(G, o)}$ denote the conditional law of $H$ given $(G, o)$.
Assume that $\alpha := \Psubbig {(\gh, \bp)}{o \in \verts(H)} > 0$ is a
constant $\bP$-a.s.
If\/ $\bP$-a.s.\
whenever $e \in \edges(\gh)$ is adjacent to $\bp$,
\rlabel e.Hsmaller
{
\Esubbig {(\gh, \bp)}{w_H(e) \bigm| e \in \edges(H)} 
\Psubbig {(\gh, \bp)}{e \in \edges(H) \bigm| \bp\in \verts(H)} 
\le
w_G(e)
\,,
}
then 
\[
\forall t > 0 \quad \Ebig{p_t(o; G)} \le \Ebig{p_t(o; H) \bigm| o \in \verts(H)}
\,.
\]
\endprocl

The case where $G$ is finite is Theorem 4.2 of \rref b.Lyons:treedom/,
although it is disguised.
The case where $G = H$ and $w_G \ge w_H$ is Theorem 5.1 of \rref b.AL:urn/.

\bigskip

\procl r.sharp
\rref t.avgrtn/ is sharp in a sense: if 
\rlabel e.wrongway
{\EBig{\sum_{e \sim o}
\Esubbig {(\gh, \bp)}{w_H(e) \bigm| e \in \edges(H)} 
\Psubbig {(\gh, \bp)}{e \in \edges(H) \bigm| \bp\in \verts(H)} }
>
\EBig{\sum_{e \sim o}
w_G(e)}
\,,
}
then for all sufficiently small, positive $t$,
\rlabel e.wrongwaysmall
{\Ebig{p_t(o; G)} > \Ebig{p_t(o; H) \bigm| o \in \verts(H)}
\,.
}
\endprocl

\bigskip
For example, let $(\gh, \bp)$ be any unimodular random rooted graph and
consider Bernoulli($\alpha$) site percolation on $\gh$.
Let $H$ be the induced subgraph. Then 
\[
\all {t \ge 0} \quad \Ebig{p_t(o; \gh)}
\le
\Ebig{p_{t/\alpha}(o; H) \bigm| o \in \verts(H)}
\,.
\]
This is obtained by using $w_H \equiv 1/\alpha$.
This is sharp: for all $\beta < \alpha$, there is some $t$ such that
$\Ebig{p_t(o; \gh)}
>
\Ebig{p_{t/\beta}(o; H)}
$.

The following corollary is immediate from \rref t.avgrtn/.

\procl c.transitive
Let $\gh$ be a unimodular transitive graph and $H$ be a random subgraph of
$\gh$ with edge weights $w_H$ such that the law of
$(H, w_H)$ is $\Aut(\gh)$-invariant. If 
\[
\all {e \sim \bp} \quad 
\Ebig{w_H(e) \bigm| e \in H} \bP[e \in H \mid o \in H]
\le
1
\,,
\]
then continuous-time simple random walk on $\gh$ and the continuous-time
network random walk on $(H, w_H)$ satisfy 
\[
\forall t > 0 \quad p_t(\bp; \gh)
\le
\Ebig{p_t(\bp; H, w_H) \bigm| o \in \verts(H)} 
\,.
\]
\endprocl

One might expect also the following as a corollary:
Suppose that $G$ is a fixed Cayley graph and $w_1$, $w_2$ are two random
fields of
positive weights on its edges with the properties that 
each field $w_i$ has an invariant law and a.s.\ $w_1(e) \ge w_2(e)$
for each edge $e$.
Then $\Ebig{p_{1, t}(o; G)} \le \Ebig{p_{2, t}(o, G)}$ for all $t
> 0$, where $p_{i, t}$ denotes the return probability to a fixed vertex $o$
at time $t$ with the weights $w_i$.
This is indeed known to be true for amenable $G$ \rref B.FontesMathieu/
and also when the pair $(w_1, w_2)$ is invariant \rref B.AL:urn/.
However, it is open in general and was asked by
Fontes and Mathieu (personal communication). 
Even more generally, the following question is open, even for finite
graphs where it was raised by \rref b.Lyons:treedom/: 

\procl q.domin
Suppose that $\bP_1$ and $\bP_2$ are two unimodular probability measures on
rooted graphs with positive edge weights for which there is a coupling
that is carried by the set of pairs $\bigl((G, o, w_G), (H, o,
w_H)\bigr)$ with $H$ a subgraph of $G$ and $w_H(e) \le w_G(e)$ for all $e
\in \edges(H)$. Is $\Esubbig 1{p_t(o; G)} \le \Esubbig 2{p_t(o; H)}$ for
all $t > 0$?
\endprocl

We prove \rref t.avgrtn/ and \rref r.sharp/ in the following section. 
Here we present the background required, especially
regarding von Neumann algebras.

We will use the notation $A \le B$ for self-adjoint operators $A$ and $B$
to mean that $B - A$ is positive semidefinite.
Sometimes we regard the edges of a graph as oriented, where we choose one
orientation (arbitrarily) for each edge. In particular, we do this whenever
we consider the $\ell^2$-space of the edge set of a graph.
In this case, we denote the tail
and the head of $e$ by $\etail e$ and $\ehead e$.
Define $d_G \colon \ell^2\big(\vertex(G)\big) \to \ell^2\big(\edges(G)\big)$ by
\[
d_G(a)(e)
:=
\sqrt{w_G(e)}\, \big[a(\etail e) - a(\ehead e)\big]
\,.
\]
Then $\Delta_G = d_G^* d_G$.

Consider the Hilbert space $\HilG := \int^\oplus \ell^2\big(\vertex(\gh)\big)
\,d\!\bP(\gh, \bp)$;
see Section 5 of \rref b.AL:urn/ for details of this direct integral.
Let $\Tr$ denote the normalized trace corresponding to $\bP$, as in
Section~5 of \rref b.AL:urn/.
That is, given an equivarant operator $T = \int^\oplus T_G \,d\!\bP(\gh,
\bp)$ on $\HilG$ in the von Neumann algebra $\alg$ associated by \rref
b.AL:urn/ to $\bP$, we define
\[
\Tr(T) := \int \ip{T_G \II{o}, \II{o}} \,d\!\bP(\gh, \bp)
\,.
\]
This trace on $\alg$ is obviously finite.
A closed densely defined operator is \dfn{affiliated}
with $\alg$ if it commutes with all unitary operators
that commute with $\alg$.
Write $\affalg$ for the set of all such operators.
An operator $T \in \affalg$ is called
\dfn{$\Tr$-measurable} if for all $\epsilon > 0$,
there is an orthogonal projection $E \in \alg$ whose
image lies in the domain 
of $T$ and $\Tr(E^\perp) < \epsilon$.
For example, $\Delta_G$ is $\Tr$-measurable because 
if $E_n$ denotes the orthogonal projection to the space of functions that
are nonzero only on those $(\gh, \bp)$
where the sum of the edge weights at $\bp$ is at most $n$,
then $\lim_{n \to\infty} \Tr(E_n^\perp) = 0$ and $\|\Delta_G
E_n\| \le 2n$.
We will need another representation of the trace.
For $s \in [0, 1]$ and a $\Tr$-measurable
operator $T \ge 0$ with spectral resolution $E_T$, define 
\[
m_s(T)
:=
\inf \bigl\{ \lambda \ge 0 \st \Tr \bigl(E_T(\lambda, \infty)\bigr) \le 1 -
s\bigr\}
\,;
\]
see Remark 2.3.1 of \rref b.FackKos/.
By Lemma 2.5(iii) of \rref b.FackKos/, if $0 \le S \le T$ are
$\Tr$-measurable, then 
\rlabel e.monotone-snum
{\all {s \in [0, 1]} \quad m_s(S) \le m_s(T)
\,.
}
A proof similar to that of Corollary 2.8 of \rref b.FackKos/ shows
that for bounded monotone $f \colon \R \to \R$ and $T \in \affalg$, we have 
\rlabel e.snum-trace
{\Tr\bigl(f(T)\bigr)
=
\int_0^1 f\bigl(m_s(T)\bigr) \,ds
\,.
}
From \rref e.snum-trace/ and \rref e.monotone-snum/, we obtain
\rlabel e.monotone-trace
{\Tr f(S) \le \Tr f(T)
}
for bounded increasing $f \colon \R \to \R$ and $0 \le S \le T$ that are
$\Tr$-measurable operators in $\affalg$. 
Furthermore, if $f$ is strictly increasing, then equality holds in \rref
e.monotone-trace/ iff $S = T$: if equality holds, then $f(S) = f(T)$
(because $\Tr$ is faithful by Lemma 2.3 of \rref b.AL:urn/),
whence $f^{-1}\bigl(f(S)\bigr) = f^{-1}\bigl(f(T)\bigr)$.

Let $w_{G, n}$ denote the weights on $G$ when for every $x \sim y$,
the edge weight $w_G(x, y)$ is replaced by 0 if
the sum of the weights incident to $x$ and $y$ is larger than $n$.
We claim that
\rlabel e.truncate
{\int p_t(\bp; w_{G, n}) \,d\!\bP(\gh, \bp)
= \lim_{n \to\infty} \int p_t(\bp; w_G) \,d\!\bP(\gh, \bp)
\,.
}
To see this,
let $E_n$ denote, as before,
the orthogonal projection to the space of functions that
are nonzero only on those $(\gh, \bp)$
where the sum of the edge weights at $\bp$ is at most $n$.
Then $\Delta_{G, w_{G, n}} E_n = \Delta_{G, w_G} E_n$ for all $n$.
Since $\lim_{n \to\infty} \Tr(E_n^\perp) = 0$, it follows that $\lim_{n
\to\infty} \Delta_{G, w_{G, n}} = \Delta_{G, w_G}$ in the measure topology
(Definition 1.5 of \rref b.FackKos/).
Since $\Delta_{G, w_{G, n}} \le \Delta_{G, w_G}$, we have $m_s\bigl(\Delta_{G, w_{G, n}}\bigr) \le
m_s\bigl(\Delta_{G, w_G}\bigr)$ by \rref e.monotone-snum/.
Therefore, $\lim_{n \to\infty} m_s\bigl(\Delta_{G, w_{G, n}}\bigr) =
m_s\bigl(\Delta_{G, w_G}\bigr)$ by Lemma 3.4(ii) of \rref b.FackKos/.
Now use $f(\lambda) := e^{-t\lambda}$ in \rref e.snum-trace/ to obtain $\lim_{n
\to\infty} \Tr\bigl(e^{-t \Delta_{G, w_{G, n}}}\bigr) = 
\Tr\bigl(e^{-t \Delta_{G, w_G}}\bigr)$, which is the same as \rref
e.truncate/.

Suppose that $\Phi$
is a positive, unital, linear map from a unital $C^*$-algebra $\algA$ to a von
Neumann algebra with finite trace, $\Tr$.
The proof of Theorem 3.9 of \rref b.Ant:Jensen/ shows that
\rlabel e.jensen
{\Tr j\big(\Phi(T)\big)
\le
\Tr \Phi\big(j(T)\big)
}
for self-adjoint operators $T \in \algA$ and functions $j
\colon \R \to \R$ that
are convex on the convex hull of the spectrum of $T$.
(In fact, those authors show the more general inequality 
$\Tr k\big(j\big(\Phi(T)\big)\big)
\le
\Tr k\big(\Phi\big(j(T)\big)\big)$ 
for every increasing convex $k$.)


\bsection{Proofs}{s.proof}

\rproof[Proof of \rref t.avgrtn/]
Suppose first that the entries of $\Delta_G$ and $\Delta_H$ are uniformly
bounded, so that $\Delta_\HilG$ and $\Delta_\HilH$ are bounded oeprators in
$\alg$.

In addition to 
the Hilbert space $\HilG := \int^\oplus \ell^2\big(\vertex(\gh)\big)
\,d\!\bP(\gh, \bp)$ we considered in the preceding section, also let
\[
\HilH := \int^\oplus 
\int^\oplus \ell^2\big(\vertex(H)\big) \,d\!\Psub {(G, o)}(H)
\,d\!\bP(\gh, \bp)
\,.
\]
By Lemma 2.3 of \rref b.AL:urn/, we have that 
\rlabel e.allalpha
{\PBig{\all {x \in
\verts(\gh)} \Psubbig {(\gh, o)}{x \in \verts(H)} = \alpha} = 1
\,.
}
Similarly, \rref e.Hsmaller/ implies that a.s.\ 
\rlabel e.Hsmaller-all
{
\all {e \in \edges(\gh)}\quad
\alpha^{-1} \int_{e \in \edges(H)} w_H(e) \,d\!\Psub {(\gh, \bp)} 
\le
w_G(e) 
\,.
}
By \rref e.allalpha/,
for every $f = \int^\oplus f(G, o) \,d\!\bP(G, o) \in \HilG$, we have that 
\[
\phi(f) := 
\alpha^{-1/2} \int^\oplus 
\int^\oplus \sum_{x \in \verts(H)} f(G, o)(x) 
\II{x}\,d\!\Psub {(G, o)}(H) \,d\!\bP(G, o)
\in \HilH
\]
has the same norm as $f$. Moreover, $\phi\colon \HilG \to \HilH$
defines an isometry, i.e., $\phi^* \phi$ is the identity map.
Define $\Phi \colon \LL(\HilH) \to \LL(\HilG)$ by 
$
\Phi T
:=
\phi^* T  \phi
$.
Then $\Phi$ is a positive unital map.

Consider quadruples $(G, H, w_G, w_H)$ of graphs $G$ and $H$ and weight
functions $w_G$ and $w_H$ with $H$ a subgraph of $G$. An isomorphism of a
pair of such quadruples is defined in the obvious way. 
As before, however, we will generally omit including the weight functions
in the notations for networks.
Similarly to how $\alg$ is defined,
let $\algA$ be the von Neumann algebra of (equivalence
classes of) bounded linear maps $T = \int^\oplus 
\int^\oplus T_{(G, o, H)} \,d\!\Psub {(G, o)}(H) \,d\!\bP(G, o)
\in \LL(\HilH)$ that are equivariant in
the sense that for all isomorphisms $\psi \colon (\gh_1, H_1)
\to (\gh_2, H_2)$, all
$\bp_1 \in \verts(\gh_1)$, $\bp_2 \in \verts(\gh_2)$, and all
$x, y \in \verts(H_1)$,
we have $(T_{(\gh_1, \bp_1, H_1)} x, y) = (T_{(\gh_2, \bp_2, H_2)} \psi x, \psi
y)$; in particular, $T_{(G, o, H)}$ does not depend on $o$.
Then $\Phi$ maps $\algA$ into $\alg$.

Let $\Delta_\HilG := \dint \Delta_G \,d\!\bP(G, o) \in \alg$ and
$\Delta_\HilH := \dint \dint \Delta_H \,d\!\Psub {(G, o)}(H) \,d\!\bP(G, o)
\in \algA$.
Then $\Phi(\Delta_\HilH) \in \alg$ and, therefore, 
$j(\Delta_\HilG), j\bigl(\Phi(\Delta_\HilH)\bigr) \in \alg$ for all bounded
Borel $j \colon \R \to \R$.

We claim that 
\rlabel e.bigger
{
\Delta_\HilG 
\ge 
\Phi(\Delta_\HilH)
\,.
}
To see this, let $f \in \HilG$.
We have 
\rlabel e.1
{
\bigpip{\Delta_\HilG(f), f}
=
\Ebig{ \|d_G f(G, o)\|^2 }
}
and 
\rlabel e.2
{
\bigpip{\Phi \Delta_\HilH(f), f}
=
\bigpip{\phi^* \Delta_\HilH \phi f, f}
=
\bigpip{\Delta_\HilH \phi f, \phi f}
\,.
}
Now
\eqaln{
\bigl(\Delta_\HilH &\phi f, \phi f\bigr)
=
\alpha^{-1} \int \!\!\! \int \|d_H f(\gh, \bp)\|^2
\,d\!\Psub {(\gh, \bp)} \,d\!\bP(\gh, \bp)
\cr&=
\alpha^{-1} \int \!\!\! \int \sum_{e \in \edges(H)} w_H(e) \bigl(f(\gh,
\bp)(\etail e) - f(\gh, \bp)(\ehead e)\bigr)^2
\,d\!\Psub {(\gh, \bp)} \,d\!\bP(\gh, \bp)
\cr&=
\alpha^{-1} \int \sum_{e \in \edges(G)} \int_{\edges(H) \ni e} w_H(e)
\,d\!\Psub {(\gh, \bp)} \cdot \bigl(f(\gh, \bp)(\etail e) - f(\gh, \bp)(\ehead e)\bigr)^2
\,d\!\bP(\gh, \bp)
\cr&\le
\int \sum_{e \in \edges(G)} w_G(e) \cdot
\bigl(f(\gh, \bp)(\etail e) - f(\gh, \bp)(\ehead e)\bigr)^2
\,d\!\bP(\gh, \bp)
\cr&=
\Ebig{ \|d_G f(G, o)\|^2 }
}
by \rref e.Hsmaller-all/.
Combining this with \rref e.1/ and \rref e.2/, we get our claimed inequality
\rref e.bigger/.

By \rref e.bigger/ and \rref e.monotone-trace/, we have
\[
\Tr j(\Delta_\HilG)
\le
\Tr j\bigl(\Phi(\Delta_\HilH)\bigr)
\]
for every decreasing function $j$.
(We have strict inequality if $j$ is strictly decreasing and we have strict
inequality in \rref e.bigger/.)
Use $j(s) := e^{-t s}$ in this and in \rref e.jensen/ to obtain 
\rlabel e.done
{
\Tr j(\Delta_\HilG)
\le
\Tr \Phi\big(j(\Delta_\HilH)\big)
\,.
}
The left-hand side equals 
$
\Ebig{p_t(o; G)} 
$.
The right-hand side equals
\eqaln{
\Tr \Phi\big(j(\Delta_\HilH)\big)
&=
\alpha^{-1} \int \!\!\! \int_{\bp \in \verts(H)} \bigip{j(\Delta_H) \II{\bp},
\II{\bp}} \,d\!\Psub {(\gh, \bp)} \,d\!\bP(\gh, \bp)
\cr&=
\Ebig{p_t(o; H) \bigm| o \in \verts(H)}
\,,}
which completes the proof of the theorem in the case of bounded vertex
weights.

We deduce the general case from this by a truncation argument.
Recall \rref e.truncate/ and its notation, which we use also for $H$.
Let $\mu_n$ be the law of $(w_{G, n}, w_{H, n})$. Since the diagonal
entries of $\Delta_{(G, w_{G, n})}$ and $\Delta_{(H, w_{H, n})}$ are
bounded and \rref e.Hsmaller/ holds $\mu_n$-a.s., we have proved that 
\[
\forall t > 0 \quad \Ebig{p_t\bigl(o; (G, w_{G, n})\bigr)} \le
\Ebig{p_t\bigl(o; (H, w_{H, n})\bigr) \bigm| o \in \verts(H)}
\,.
\]
Taking $n \to\infty$ and using the bounded convergence theorem,
we get the desired result.
\Qed

A similar proof shows that \rref e.done/ holds
if $j$ is any decreasing convex function.

\rproof[Proof of \rref r.sharp/]
The right-hand side of \rref e.wrongway/ is equal to
$\Ebig{\Delta_G(o, o)}$ and the left-hand side is
$\Ebig{\Delta_H(o, o) \bigm| o \in \verts(H)}$.
Now both sides of \rref e.wrongwaysmall/
equal 1 for $t = 0$. We claim that the derivative of
the left-hand side at $t = 0$ is
$-\Ebig{\Delta_H(o, o) \bigm| o \in \verts(H)}$
and the derivative of the right-hand side at $t = 0$ is
$-\Ebig{\Delta_G(o, o)}$.
This clearly implies the remark.
To evaluate these derivatives, note that
for every fixed $G$, the spectral representation 
\[
p_t(o; G)
=
\int_0^\infty e^{-\lambda t} \,d\bigip{E_{\Delta_G}(\lambda) \II{o},
\II{o}}
\]
shows that $t \mapsto p_t(o; G)$ is monotone decreasing and convex.
By Tonelli's theorem, it follows that for $(G, o) \sim \bP$,
\[
\Ebig{p_t(o; G)} - 1
=
\EBig{\int_0^t p'_s(o; G) \,ds}
=
\int_0^t \Ebig{p'_s(o; G)} \,ds
\,.
\]
The fundamental theorem of calculus and the monotone convergence theorem
now yield that 
\[
\frac{d}{dt} \Ebig{p_t(o; G)} \biggr|_{t=0}
=
\Ebig{p'_0(o; G)}
=
-\Ebig{\Delta_G(o, o)}
\,.
\]
A similar calculation applied to the distribution of $(H, o)$ given $o \in
\verts(H)$ yields the derivative of the left-hand side of \rref
e.wrongwaysmall/.
\Qed

\bibliographystyle{amsalpha}
\bibliography{\jobname.bib}

\end{document}